\documentclass[12pt]{article}
\usepackage[fleqn]{amsmath}
\usepackage{latexsym}
\usepackage[dvips]{graphicx}
\usepackage{psfrag}
\newtheorem{thm}{Theorem}
\newtheorem{lem}[thm]{Lemma}

\newtheorem{con}[thm]{Conjecture}

\def\qed{$\hfill \Box$}
\def\prf{\emph{Proof: }} 
\def\b{\backslash}
\def\con{/}

\def\indr{\mc{I}^{(r)}(G)}

\def\del{\!\downarrow\!}
\def\tr{\textrm}

\def\A{\mathcal{A}}
\def\B{\mathcal{B}}
\def\C{\mathcal{C}}
\def\D{\mathcal{D}}
\def\E{\mathcal{E}}
\def\F{\mathcal{F}}
\def\I{\mathcal{I}}
\def\mc{\mathcal}
\begin{document}
\title{Compression and Erd\H os-Ko-Rado graphs}
\author{Fred Holroyd$^{1}$ and John Talbot$^{2}$\\[3mm]
$^{1}$Department of Pure Mathematics, The Open University,\\
Walton Hall, Milton Keynes MK7~6AA, United Kingdom\\
\texttt{f.c.holroyd@open.ac.uk}\\ \\
$^{2}$Merton College and the Mathematical Institute,\\
Oxford University, OX1 4JD\\
\texttt{talbot@maths.ox.ac.uk}}
\maketitle
\begin{abstract}
For a graph $G$ and integer $r\geq 1$ we denote the collection of independent $r$-sets of $G$ by $\I^{(r)}(G)$. If $v\in V(G)$ then $\I_v^{(r)}(G)$ is the collection of all independent $r$-sets containing $v$. A graph $G$, is said to be \emph{$r$-EKR}, for $r\geq 1$, iff no intersecting family $\A\subseteq \I^{(r)}(G)$ is larger than $\max_{v\in V(G)}|\I^{(r)}_v(G)|$. There are various graphs which are known to have his property: the empty graph of order $n\geq 2r$ (this is the celebrated Erd\H os-Ko-Rado theorem), any disjoint union of at least $r$ copies of $K_t$ for $t\geq 2$, and any cycle. In this paper we show how these results can be extended to other classes of graphs via a compression proof technique.

In particular we show that any disjoint union of at least $r$ complete graphs, each of order at least two, is $r$-EKR. We also show that paths are $r$-EKR for all $r\geq 1$.
\end{abstract}
\section{Introduction}
An \emph{independent} set in a graph $G=(V,E)$, is a subset of the vertices not containing any edges. For an integer $r\geq 1$ we denote the collection of independent $r$-sets of $G$ by\[
\I^{(r)}(G)=\{A\subset V(G):|A|=r\tr{ and $A$ is an independent set}\}.\]
If $v\in V(G)$ then the collection of independent $r$-sets containing $v$ is \[
\I_v^{(r)}(G)=\{A\in \I^{(r)}(G):v\in A\}.\] Such a family is called a \emph{star}.

A graph $G$ is \emph{$r$-EKR} iff no intersecting family of independent $r$-sets is larger than the largest star in $\I^{(r)}(G)$. If $G$ is $r$-EKR and any intersecting family $\A\subseteq\I^{(r)}(G)$ of maximum size is a star then $G$ is said to be \emph{strictly $r$-EKR}.

In this setting the classical Erd\H os-Ko-Rado theorem can be stated as follows.
\begin{thm}[Erd\H os-Ko-Rado\cite{EKR}]\label{EKR}If $E_n$ is the empty graph of order $n$ then $E_n$ is $r$-EKR for $n\geq 2r$ and strictly so for $n>2r$. 
\end{thm}
There are a number of other recent results of this type. 
\begin{thm}[Bollob\'as and Leader\cite{BL}]\label{BLthm}
If $n \geq r$, $t\geq 2$ and $G$ is the disjoint union of $n$ copies of $K_t$ then $G$ is $r$-EKR and strictly so unless $t=2$ and $n=r$.
\end{thm}
The following result was previously a conjecture of Holroyd and Johnson \cite{HJ}.
\begin{thm}[Talbot\cite{JT}]\label{cyclethm}
For $1\leq k \leq n$ the $k^{\tr{th}}$ power of the $n$-cycle, $C_n^k$, is the graph with vertex set $[n]=\{1,2,\ldots,n\}$ and edges between $a,b \in [n]$ iff $1\leq |a-b\mod n|\leq k$. Then $C_n^k$ is $r$-EKR for all $r\geq 1$ and strictly so unless $n=2r+2$ and $k=1$.
\end{thm}

Although this is not made explicit in \cite{JT}, the proof of Theorem \ref{cyclethm} uses a type of compression that is essentially equivalent to contracting an edge in the underlying graph. In the present paper we wish to show how this idea can be used to show that various other graphs are also $r$-EKR.

In particular we have the following result extending Theorem \ref{BLthm}.

\begin{thm}\label{completethm}If $G$ is the disjoint union of $n\geq r$ complete graphs each of order at least two then $G$ is $r$-EKR.
\end{thm} 

We also show that an analogue of Theorem \ref{cyclethm} holds for paths.
\begin{thm}\label{paththm}
For $1\leq k \leq n$ the $k^{\tr{th}}$ power of the $n$-path, $P_n^k$, is the graph with vertex set $[n]=\{1,2,\ldots,n\}$ and edges between $a,b \in [n]$ iff $1\leq |a-b|\leq k$. Then $P_n^k$ is $r$-EKR for all $r\geq 1$.
\end{thm}

Our main results are Theorems \ref{completethm} and \ref{paththm} but the proof technique also extends to other types of graph and our final theorem gives
 a large class of graphs which are all $r$-EKR.

\begin{thm}\label{messthm}
If $G$ is a disjoint union of $n\geq 2r$ complete graphs, cycles and paths, including an isolated singleton, then $G$ is $r$-EKR.\end{thm}

In  order to state the two key lemmas we require some notation.

If $e$ is an edge of the graph $G=(V,E)$ we define $G/e$ to be the graph obtained from $G$ by contracting the edge $e$. We also define $G\del e$ to be the graph obtained from $G$ by removing the vertices in $e$ as well as their neighbours. As usual we denote the neighbours of a vertex $v$ by $\Gamma(v)$.

The following two technical lemmas relate intersecting families and stars in $\I^{(r)}(G)$ to intersecting families and stars in $\I^{(r)}(G/e)$ and $\I^{(r-1)}(G\del e)$. These will enable us to prove our main results by induction.

\begin{lem}
\label{partlem}
Let $G=(V,E)$ be a graph and $\A\subseteq \mc{I}^{(r)}(G)$ be an intersecting family. If $e=\{v,w\}\in E$ is an edge in $G$ then there exist families $\B,\C,\D$ and $\E$ satisfying:
\begin{itemize}
\item[(i)] $|\A|=|\B|+|\C|+|\D|+|\E|$.
\item[(ii)] $\B\subseteq \mc{I}^{(r)}(G\con e)$ is intersecting.
\item[(iii)] $\C\subseteq \mc{I}^{(r-1)}(G\del e)$ is intersecting.
\item[(iv)] $\D=\{A\in\A:v\in A\tr{ and }\Gamma(w)\cap (A\b\{v\})\neq \emptyset\}$.
\item[(v)] $\E=\{A\in\A:w\in A\tr{ and }\Gamma(v)\cap (A\b\{w\})\neq \emptyset\}$.
\item[(vi)] If $C\in \C$ and $F\in\D\cup \E$ then $C\cap F\cap V(G\del e)\neq \emptyset$.
\item[(vii)] If $D\in \D$ and $E\in\E$ then $D\cap E\cap V(G\del e)\neq \emptyset$.
\end{itemize}
\end{lem}

\begin{lem}
\label{starlem}If $e=\{v,w\}$ is an edge in the graph $G=(V,E)$ and $x\in V(G\del e)$ then
\[
|\mc{I}_x^{(r)}(G)|=|\mc{I}_x^{(r)}(G\con e)|+ |\mc{I}_x^{(r-1)}(G\del e)|+|\D_x|+|\E_x|.\] Where \[\D_x=\{A\in \mc{I}_x^{(r)}(G):v\in A\tr{ and }\Gamma(w)\cap (A\b\{v\})\neq \emptyset\}\] and \[\E_x=\{A\in \mc{I}_x^{(r)}(G):w\in A\tr{ and }\Gamma(v)\cap (A\b\{w\})\neq \emptyset\}.\]
\end{lem}

\section{Proofs}
\emph{Proof of Lemma \ref{partlem}:}
Let $\A\subseteq \indr$ be intersecting. We consider the effect of contracting an edge $e=\{v,w\}$ in $G$. We define a contraction function, $c:V(G)\to V(G/e)$ by \[
c(x)=\left\{\begin{array}{ll}v, & x=w,\\ x,& x\neq w.\end{array}\right.\]
Let \begin{eqnarray*}
\B&=&\{c(A):A\in \A\tr{ and $c(A)\in\I^{(r)}(G/e)$}\},\\
 \C&=&\{A\b\{v\}:v\in A\in\A\tr{ and }A\b\{v\}\cup\{w\}\in\A\},\\
\D&=&\{A\in\A:v\in A\tr{ and }\Gamma(w)\cap (A\b\{v\})\neq \emptyset\},\\
\E&=&\{A\in\A:w\in A\tr{ and }\Gamma(v)\cap (A\b\{w\})\neq \emptyset\}.
\end{eqnarray*}
If $A,B\in\A$ then $c(A)=c(B)$ iff $A\Delta B=\{v,w\}$. Hence\[ |\{A\in\A:c(A)\in\I^{(r)}(G/e)\}|=|\B|+|\C|.\] Also if $A\in\A$ then $c(A)\not\in \I^{(r)}(G/e)$ iff $A\in \D\cup \E$. Hence $|\A|=|\B|+|\C|+|\D|+|\E|$, which is (i).

The fact that $\B\subseteq\I^{(r)}(G/e)$ is  an intersecting family follows simply because $\A$ is intersecting, so (ii) holds.  

If $C\in\C$ then $C\cup\{v\},C\cup\{w\}\in \A$ hence $C\in I^{(r-1)}(G\del e)$.
With a little more thought it is also clear that $\C$ is an intersecting family. Let $C,D\in\C$, if $C\cap D=\emptyset$ then $\A$ contains the two disjoint sets $C\cup\{v\}$ and $D\cup\{w\}$. This contradicts the fact that $\A$ is intersecting. Hence $\C$ is also intersecting, and so (iii) holds.

The definitions of the families $\D$ and $\E$ give (iv) and (v).

To see that (vi) holds let $C\in \C$, so $C\cup\{v\},C\cup\{w\}\in \A$. If $F\in\D\cup\E\subseteq \A$ then $(C\cup\{w\})\cap F\neq\emptyset$ and $(C\cup\{v\})\cap F\neq\emptyset$ but either $v\not\in F$ or $w \not\in F$. Hence $C\cap F\cap V(G\del e)\neq\emptyset$.

Finally if $D\in\D$ and $E\in \E$ then $v\in D$ and $w\in E$ imply that \[D\cap E\cap (\Gamma(v)\cup\Gamma(w)\cup\{v,w\})=\emptyset.\] So $D\cap E\neq\emptyset$ implies that (vii) must hold. \qed

\emph{Proof of Lemma \ref{starlem}:} This follows similarly to the proof of Lemma \ref{partlem}, via contracting the edge $e=\{v,w\}$. Let $c:V(G)\to V(G/e)$ be as defined in the proof of Lemma \ref{partlem}. 

Then $c$ is a surjection between the families $\I_x^{(r)}(G)\b(\D_x\cup \E_x)$ and $\I_x^{(r)}(G/e)$. Moreover $c(A)=c(B)$ iff $A\Delta B=\{v,w\}$ and the number of sets in $\I_x^{(r)}(G/e)$ with two preimages under $c$ is exactly $|\I_x^{(r-1)}(G\del e)|$. Hence
\[
|\mc{I}_x^{(r)}(G)|=|\mc{I}_x^{(r)}(G\con e)|+ |\mc{I}_x^{(r-1)}(G\del e)|+|\D_x|+|\E_x|.\] 
\qed

\emph{Proof of Theorem \ref{completethm}:} We prove this result by induction on $r$. It is clearly true for $r=1$ so we may suppose that $r\geq 2$ and the result holds for smaller values of $r$. 

We now use induction on the number of vertices in $G$. Theorem \ref{BLthm} implies that the result holds when $G$ consists of $n\geq r$ copies of $K_t$, for $t\geq 2$. So let \[G=K_{t_1}\cup\cdots\cup K_{t_n},\] with $2\leq t_1\leq t_2\leq \cdots \leq t_n$, not all equal. We may suppose that the result holds for all graphs of the correct form with fewer vertices.

Suppose that $\A\subseteq\I^{(r)}(G)$ is intersecting. Let $v,w\in K_{t_n}$, we will  contract the edge $e=\{v,w\}$. Then \[
G/e= K_{t_1}\cup\cdots\cup K_{t_n-1},\] and \[ G\del e= K_{t_1}\cup\cdots\cup K_{t_{n-1}}.\] 
Using the notation of Lemma \ref{partlem} we have $\D=\E=\emptyset$. Hence by Lemma \ref{partlem} (i) 
\begin{equation}\label{eq1}
|\A|=|\B|+|\C|.\end{equation}  Then for any $x\in K_{t_1}\subseteq G\del e$ we have, in notation of Lemma \ref{starlem}, that $\D_x=\E_x=\emptyset$. So Lemma \ref{starlem} implies that \begin{equation}\label{eq2}|\I_x^{(r)}(G)|=|\I_x^{(r)}(G/e)|+|\I_x^{(r-1)}(G\del e)|.\end{equation}

The observation that $t_1\leq t_i$, for any $1\leq i\leq n$, implies that we also have \begin{equation}\label{eq3}
 |\I_x^{(r)}(G/e)|=\max_{v\in V(G/e)}|\I_v^{(r)}(G/e)|\end{equation} and
\begin{equation}\label{eq4}
 |\I_x^{(r-1)}(G\del e)|=\max_{v\in V(G\downarrow e)}|\I_v^{(r-1)}(G\del e)|.\end{equation}
Now $t_n\geq 3$ so $G/e$ is a smaller graph of the correct form and hence is $r$-EKR. Then Lemma \ref{partlem} (ii) and (\ref{eq3}) imply that \begin{equation}\label{eq5}
|\B|\leq |\I_x^{(r)}(G/e)|.\end{equation} Also $G\del e$ is $(r-1)$-EKR, since the result holds for smaller values of $r$. So Lemma \ref{partlem} (iii) and (\ref{eq4}) imply that \begin{equation}\label{eq6}
|\C|\leq |\I_x^{(r-1)}(G\del e)|\end{equation} 
Hence, using equations (\ref{eq1}), (\ref{eq2}), (\ref{eq5}) and (\ref{eq6}), we obtain
\[
|\A|\leq |\I_x^{(r)}(G)|.\]\qed

\emph{Proof of Theorem \ref{paththm}:}
We first note that for any $n$, $r$ and $k$\[\max_{x\in V(P_n^k)}|\mc{I}_x^{(r)}(P_n^k)|\] is achieved  by taking $x\in\{1,n\}$. 

Again we prove this result by induction on $r$. The result clearly holds for $r=1$ so we may assume $r\geq 2$ and that the result is true for smaller values of $r$. 

We now prove the result for $r$ by induction on $n$. For $n<(r-1)k+r$ there is nothing to prove since $\I^{(r)}(P_n^k)$ is empty. For $n=(r-1)k+r$ the result also holds (since there is only one set in $\I^{(r)}(P_n^k)$). So we may assume that $n\geq (r-1)k+r+1$ and that the result holds for smaller values of $n$. In particular $n\geq k+3$.

Let $\A\subseteq \mc{I}^{(r)}(P_n^k)$ be intersecting. Set $w=n$, $v=n-1$ and $e=\{n-1,n\}$, and apply Lemma \ref{partlem}. Let $\B,\C,\D$ and $\E$ be the families defined in Lemma \ref{partlem}. In this case $G/e$ is $P_{n-1}^k$, while $G\del e$ is $P_{n-k-2}^k$. Note that $n\geq k+3$ implies that $G\del e$ is non-empty.

We see that in this case $\D$ is empty and \[\E=\{A\in\A: n,n-k-1\in A\}.\] Let \[\F=\{A\b\{n\}:A\in\E\}\] and consider $\C\cup\F$. Note that this is a disjoint union since $n-k-1$ belongs to every set in $\F$ but to no set in $\C$. Hence \begin{equation}\label{eq6b}
|\C\cup \F|=|\C|+|\F|=|\C|+|\E|.\end{equation}

Parts (iii) and (vi) of Lemma \ref{partlem} imply that $\C\cup\F$ is an intersecting family of independent $(r-1)$-sets in the subgraph of $P_n^k$ induced by $\{1,2,\ldots,n-k-1\}$, which is $P_{n-k-1}^k$. Our inductive hypothesis for $r$ then implies that \begin{equation}\label{eq7}|\C\cup\F|\leq |\mc{I}_1^{(r-1)}(P_{n-k-1}^k)|.\end{equation}
Now $G/e$ is $P_{n-1}^k$, so part (ii) of Lemma \ref{partlem} and our inductive hypothesis for $n$ imply that \begin{equation}\label{eq8}|\B|\leq |\mc{I}_1^{(r)}(P_{n-1}^k)|.\end{equation}
Lemma \ref{partlem} (i), together with equations (\ref{eq6b}), (\ref{eq7}) and (\ref{eq8}) imply that
\begin{equation}\label{eq9}
|\A|=|\B|+|\C|+|\D|+|\E|\leq |\I_1^{(r)}(P_{n-1}^k)|+|\I_1^{(r-1)}(P_{n-k-1}^{k})|.\end{equation}

Applying Lemma \ref{starlem} we obtain \begin{equation}\label{eq10}|\mc{I}_1^{(r)}(P_n^k)|=|\mc{I}_1^{(r)}(P_{n-1}^k)|+|\mc{I}_1^{(r-1)}(P_{n-k-2}^k)|+|\E_1|,\end{equation}
where \[
\E_1=\{A\in \I_1^{(r)}(P_n^k):n-k-1,n\in A\}.\]
Then it is easy to check that 
\[
|\I_1^{(r-1)}(P_{n-k-1}^k)|=|\I_1^{(r-1)}(P_{n-k-2}^k)|+|\E_1|,\]

  and so equations (\ref{eq9}) and (\ref{eq10}) imply that \[|\A|\leq|\mc{I}_1^{(r)}(P_n^k)|.\] as required.
\qed

We turn finally to a proof of Theorem \ref{messthm}. The key ideas have already been presented in Lemmas \ref{partlem} and \ref{starlem} as well as in the proofs of Theorems \ref{completethm} and \ref{paththm}. For this reason our proof is essentially a sketch. %The interested reader should easily be able to fill in the details.

\emph{Proof of Theorem \ref{messthm}:} We will say that a graph $G$ is \emph{$r$-mixed}, for an integer $r\geq 1$, if it satisfies the conditions of Theorem \ref{messthm}. So $G$ is $r$-mixed iff it is the disjoint union of at least $2r$ complete graphs, paths and cycles, including at least one isolated singleton.

We prove the result by induction on $r$. It is clearly true for $r=1$ so we may suppose $r\geq 2$ and that the result holds for smaller values of $r$.

We now prove the result for $r$ by induction on $|V(G)|$. If $G$ is $r$-mixed then $|V(G)|\geq 2r$ with equality iff $G=E_{2r}$.  So if $|V(G)|=2r$ then Theorem \ref{EKR} implies that $G$ is $r$-EKR. Hence we may suppose that $|V(G)|> 2r$ and that any $r$-mixed graph with fewer vertices is also $r$-EKR.

Now either $G$ is an empty graph of order at least $2r+1$, in which case the result holds by Theorem \ref{EKR}, or $G$ contains an edge. So we may suppose that $G$ contains an edge $e=\{v,w\}$. We also know that $G$ contains an isolated singleton $x$. 

It is easy to check that if $H$ is one of the three graphs $G$, $G/e$ or $G\del e$,  and $s\leq r$ then \begin{equation}\label{lasteq}
|\I_x^{(s)}(H)|=\max_{v\in V(H)}|\I_v^{(s)}(H)|.\end{equation}
Also $G/e$ is $r$-mixed and has less vertices than $G$, while $G\del e$ is $(r-1)$-mixed. Hence  by our two inductive hypotheses $G/e$ is $r$-EKR and $G\del e$ is $(r-1)$-EKR. The result then follows by applying Lemmas \ref{partlem} and \ref{starlem} together with (\ref{lasteq}). The details depend on whether $e$ was an edge from a complete graph, a cycle or a path but follow similarly to the proofs of Theorems \ref{cyclethm}, \ref{completethm} and \ref{paththm}.\qed 

\end{document}